\documentclass[a4paper,12pt]{article}

\usepackage{amsfonts}
\usepackage{latexsym}
\usepackage{epsfig}
\usepackage{amssymb}
\usepackage{graphicx}
\usepackage{oldgerm}
\usepackage{theorem}
\theorembodyfont{\itshape}
\newtheorem{thm}{Theorem}[section]
\newtheorem{lem}[thm]{Lemma}

\newtheorem{prop}[thm]{Proposition}

\newcommand{\mib}[1]{\mbox{\boldmath $#1$}}
\def\B{\mib{B}}

\def\C{{\bf C}}
\def\R{{\bf R}}
\def\Z{{\bf Z}}
\def\W{{\bf W}}

\def\E{{\bf E}}
\def\P{{\bf P}}

\def\BES3{${\rm BES}_3$}

\def\X{\mib{X}}
\def\x{\mib{x}}
\def\y{\mib{y}}
\newcommand{\SSC}[1]{\section{#1}\setcounter{equation}{0}}

\newcommand{\qed}{\hbox{\rule[-2pt]{3pt}{6pt}}}

\begin{document}

\title{\bf Two Bessel Bridges Conditioned \\
Never to Collide, Double Dirichlet Series, \\
and Jacobi Theta Function}

\author{
Makoto Katori $\cdot$
Minami Izumi $\cdot$
Naoki Kobayashi
\footnote{
M. Katori $\cdot$ M. Izumi $\cdot$ N. Kobayashi:
Department of Physics,
Faculty of Science and Engineering,
Chuo University, 
Kasuga, Bunkyo-ku, Tokyo 112-8551, Japan 
}
\footnote{
M. Katori:
e-mail: katori@phys.chuo-u.ac.jp}
\footnote{
M. Izumi:
e-mail: izumi@phys.chuo-u.ac.jp}
\footnote{
N. Kobayashi:
e-mail: knaoki@phys.chuo-u.ac.jp
}}
\date{13 March 2008}
\pagestyle{plain}
\maketitle
\begin{abstract}
It is known that the moments of the maximum
value of a one-dimensional conditional Brownian motion, 
the three-dimensional Bessel bridge
with duration 1 started from the origin, 
are expressed using the Riemann zeta function.
We consider a system of two Bessel bridges,
in which noncolliding condition is imposed.
We show that the moments of the maximum value
is then expressed using the double Dirichlet
series, or using the integrals of products of the
Jacobi theta functions and its derivatives.
Since the present system will be provided as
a diffusion scaling limit of a version of
vicious walker model, 
the ensemble of 2-watermelons with
a wall, the dominant terms in long-time
asymptotics of moments of height of
2-watermelons are completely determined.
For the height of 2-watermelons with a wall, 
the average value was recently
studied by Fulmek by a method of
enumerative combinatorics. 
\end{abstract}

\SSC{Introduction}

Let $\B(t)=(B_1(t), B_2(t), B_3(t)), t \geq 0$ 
be the three-dimensional
Brownian motion (BM), in which three components
$B_j(t), j=1,2,3$ are given by independent 
one-dimensional standard BMs.
The three-dimensional Bessel process (\BES3),
$X(t)$, started from $x > 0$ is defined as the
radial part of $\B(t)$,
\begin{eqnarray}
X(t) &\equiv& |\B(t)| \nonumber\\
&=& \sqrt{B_1(t)^2+B_2(t)^2+B_3(t)^2}, \quad
t \geq 0
\nonumber
\end{eqnarray}
with $X(0)=x$.
\BES3 is a diffusion process on 
$\R_{+}=\{x \in \R : x \geq 0\}$,
where $\R$ denotes the set of all real numbers.
By It\^o's formula 
we can show that it satisfies the stochastic differential
equation of the form
$$
dX(t)=dB(t)+\frac{1}{X(t)}dt, \quad
t \geq 0, \quad X(0)=x,
$$
where $B(t)$ is the one-dimensional standard BM
different from $B_j(t)$'s used to give $\B(t)$ above.
We can prove that $X(t) \to \infty$ 
in $t \to \infty$ with probability one
for all $x \geq 0$,
{\it i.e.} \BES3 is transient.
For the basic properties of \BES3,
see, for example,
3.3 C in \cite{KS91}, VI.3 in \cite{RY98}, IV.34 in \cite{BS02}.

The {\it three-dimensional Bessel bridge 
with duration 1 started from the origin},
$\widetilde{X}(t), t \in [0,1]$, is
then defined as the \BES3 conditioned
$$
x=X(0) = 0 \quad \mbox{and} \quad
X(1)=0.
$$
\begin{figure}
\includegraphics[width=0.8\linewidth]{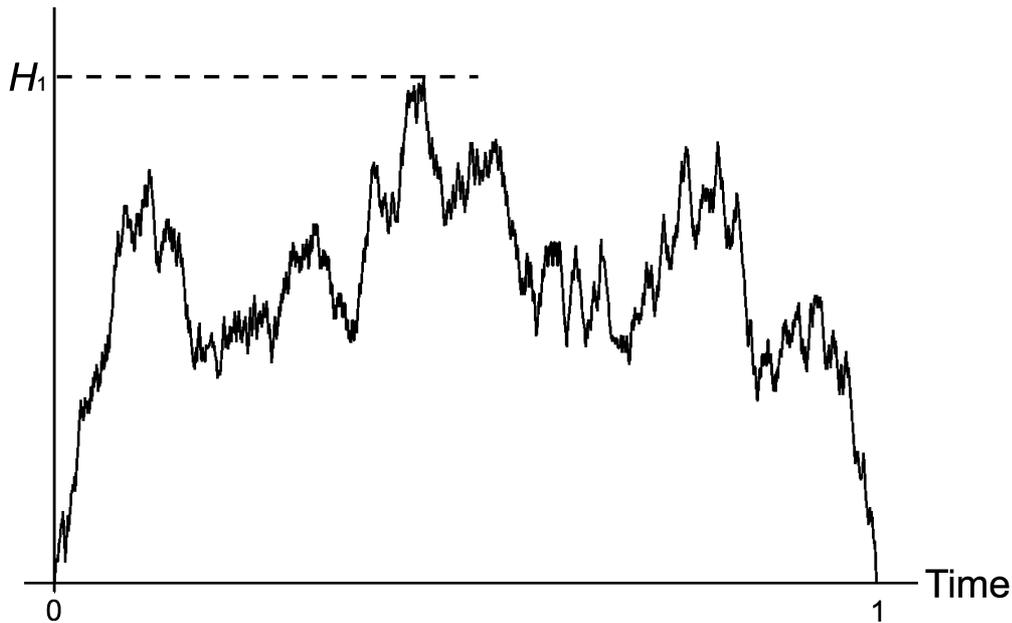}
\caption{Sample path of three-dimensional Bessel bridge with duration 1.
\label{fig:BridgeFig1}}
\end{figure}
Figure 1 illustrates a sample path
of $\widetilde{X}(t)$
on the spatio-temporal plane
$(t,x) \in [0,1] \times \R_{+}$.
In \cite{BPY01}, a variety of probability laws
associated with conditional Brownian motions
are discussed, which are related to
the Jacobi theta function and the Riemann
zeta function.
One of them is the probability law
of the maximum value of
$\widetilde{X}(t)$;
\begin{equation}
H_1 \equiv \max_{0 < t < 1} \widetilde{X}(t).
\label{eqn:H1}
\end{equation}
Let $\E[H_1^{s}]$ be the $s$-th moment of $H_1$.
The following equality is discussed in
\cite{BPY01},
\begin{equation}
\E[H_1^{s}]=2 \left(\frac{\pi}{2} \right)^{s/2}
\xi(s), \quad s \in \C,
\label{eqn:equal1}
\end{equation}
where
$\C$ denotes the set of all complex numbers, and
$$
\xi(s)=\frac{1}{2}s(s-1) \pi^{-s/2}
\Gamma(s/2) \zeta(s)
$$
with the gamma function
\begin{equation}
\Gamma(s)=\int_{0}^{\infty} du \, 
u^{s-1} e^{-u}, \quad \Re s > 0,
\label{eqn:gamma1}
\end{equation}
and with the Riemann zeta function
$$
\zeta(s)=\sum_{n=1}^{\infty} \frac{1}{n^s},
\quad \Re s > 1.
$$
See also  Chapter 11 in \cite{Yor97}.

We know the two facts;
(i) the BM can be realized as the diffusion scaling
limit of the simple random walk,
(ii) the probability law of \BES3 is
equal to that of the BM conditioned to stay
positive.
Combination of them will lead to the following.
For a fixed $n \geq 1$, 
consider one-dimensional simple random walks
started from the origin,
which visit only positive sites $\{1,2,3, \cdots\}$
up to time $2n$ and return to the origin
at time $2n$.
Sample paths of such conditional
random walks are called Dyck paths of length $n$
in combinatorics.
The height of Dyck path $h_1(2n)$ 
is defined as the maximum site
visited by the walker.
Let $\langle \, \cdot \, \rangle$ denote
the average over all Dyck paths
with uniform weight.
Then we will have the relation
\begin{equation}
\lim_{n \to \infty}
\Big\langle \left(
\frac{h_1(2n)}{\sqrt{2n}}\right)^{s} \Big\rangle
=\E[H_1^s], \quad s \in \C.
\label{eqn:equal2}
\end{equation}

The classical work of de Bruijn, Knuth and Rice
in enumerative combinatorics
\cite{dBKR72} gives
\begin{equation}
\langle h_1(2n) \rangle \simeq
\sqrt{\pi n} -\frac{3}{2} +o(1)
\quad \mbox{in} \quad n \to \infty.
\label{eqn:BKR}
\end{equation}
Here we should note that,
through the relations (\ref{eqn:equal1})
and (\ref{eqn:equal2}), 
if we only consider the dominant term 
in (\ref{eqn:BKR}) proportional to 
$\sqrt{n}$,
this result in combinatorics means nothing 
but the fact $\xi(1)=1/2$.
It is rather obvious if we know the following integral
representation of $\xi(s)$ due to Riemann,
\begin{equation}
\xi(s)=
\frac{1}{2}+\frac{1}{4}s(s-1)
\int_{1}^{\infty} du \,
(u^{s/2-1}+u^{(1-s)/2-1})
(\vartheta(u)-1),
\label{eqn:xisym1}
\end{equation}
where $\vartheta(u)$ is a version of
the Jacobi theta function
\begin{equation}
\vartheta(u)=\sum_{n=-\infty}^{\infty}
e^{-\pi n^2 u}, \quad u > 0.
\label{eqn:theta1}
\end{equation}

Recently Fulmek reported a generalization of the
result of de Bruijn, Knuth and Rice,
by calculating the asymptotics of the average
height of {\it 2-watermelons with a wall} \cite{Ful07}.
In general, the uniform ensemble of 
$N$-watermelons, $N \geq 2$, is a version of
vicious walker model of Fisher \cite{Fis84}.
In this version
the starting points and the ending points
of $N$ vicious walkers 
({\it i.e.} nonintersecting random walks) 
are fixed to the sites located near to the origin.
When we impose the condition to stay positive
for all vicious walkers, we say
``with a wall" (at the origin)
\cite{AME91,EG95,KGV00,Gillet03,Kra06,Kuij08}.
The height of $N$-watermelon is the maximum site
visited by the vicious walker, who walks
the farthest path from the origin.
Let $h_2(2n)$ be the height of 2-watermelon
with a wall. Fulmek showed
\begin{equation}
\langle h_2(2n) \rangle \simeq
c_2 \sqrt{n} -\frac{3}{2}+o(1)
\quad \mbox{in} \quad n \to \infty
\quad \mbox{with} \quad
c_2=2.57758 \cdots.
\label{eqn:Fulmek1}
\end{equation}
Here the factor $c_2=2.57758 \cdots$ of the dominant term
proportional to $\sqrt{n}$ was given by numerical 
evaluation of the ``constant terms" in 
Laurent expansions of a version of double Dirichlet
series. The terms are represented by
integrals of functions expressed using
the Jacobi theta function (\ref{eqn:theta1}) and
its derivatives.
It should be emphasized the fact that
Fulmek succeeded in proving the $N=2$ case
of the conjecture of Bonichon and Mosbah \cite{BM03},
$$
\langle h_N(2n) \rangle
\simeq \sqrt{(1.67N-0.06) 2n} \quad
\mbox{in} \quad n \to \infty
$$
obtained by computer simulations
for the average height $h_{N}(2n)$
of general $N$-watermelons with a wall,
$N \geq 1$.
It seems to be highly nontrivial to extend
his method to evaluate the asymptotics of
higher moments $\langle h_N(2n)^s \rangle,
s \geq 2$ for $N=2$ and $N \geq 3$.
See the paper by Feierl on the recent progress
in this combinatorial method \cite{Feierl07}.

Here we propose a different
method to calculate the dominant terms of
all moments of height for 2-watermelons
with a wall.
We will perform the diffusion scaling limit first.
Following the argument of \cite{KT02,KT03a,Gillet03},
we can prove that the diffusion scaling limit of
the $N$-watermelons with a wall
provides the noncolliding system of $N$ Bessel bridges,
$\widetilde{\X}=(\widetilde{X}_1(t), \cdots, 
\widetilde{X}_N(t)) 
\in \W_N^{\rm C} \equiv 
\{(x_1, \cdots, x_N): 0 < x_1 < \cdots < x_N\},
0 < t < 1$. It implies
\begin{equation}
\lim_{n \to \infty} \Big\langle
\left( \frac{h_{N}(2n)}{\sqrt{2n}} \right)^s \Big\rangle
= \E [ H_N^{s} ], \quad N \geq 2,
\label{eqn:limitA}
\end{equation}
where
\begin{equation}
H_N= \max_{0 < t < 1}
\widetilde{X}_N(t).
\label{eqn:HN}
\end{equation}
In the present paper we determine
$\E[H_2^{s}]$ for arbitrary $s$ for the 
two Bessel bridges with noncolliding condition.
Noncolliding diffusion particle systems
are interesting and important statistical-mechanical
processes,
since they are related to the group representation-theory, 
the random matrix theory, 
and the exactly solved nonequilibrium 
statistical-mechanical models
({\it e.g.}, ASEP and polynuclear growth models)
\cite{KT07b}.
The present system of noncolliding Bessel bridges
is related to the class C ensemble of
random matrices discussed by
Altland and Zirnbauer \cite{AZ96,AZ97}
(see Sect.V.C of \cite{KT04})
and it is a special case with
parameters $(\nu, \kappa)=(1/2, 3)$
of the noncolliding generalized meanders \cite{KT07a}
(see also \cite{TW07}).

As demonstrated in \cite{KT07b}, the noncolliding
diffusion processes can be regarded as the 
multivariate extensions of Bessel processes.
The present study suggests the possibility
that the connection between the conditional BMs
and the number theoretical functions
({\it e.g.}, the Jacobi theta function, the Riemann
zeta functions, and Dirichlet series) reported in
\cite{Yor97,BPY01} will be extended to
many particle and multivariate systems.

The paper is organized as follows.
In Sect.2 we give the precise description of
the systems and main results.
In Sect.3, we give proofs of our theorems
for the $N=2$ case and
show formulas, which are useful 
to perform the numerical evaluation
of moments.

\SSC{Models and Results}
\subsection{Reflection Principle and Karlin-McGregor Formula}
The transition probability density of the
one-dimensional standard BM is given by the
heat-kernel
$$
p(t,y|x) =
\frac{1}{\sqrt{2 \pi t}} \exp
\left\{ - \frac{(y-x)^2}{2t} \right\},
\quad x, y \in \R, \, t \geq 0.
$$
By the reflection principle of BM,
the transition probability density of the BM
with an absorbing wall at the origin
is given by
\begin{eqnarray}
p_{1}(t, y|x) &=& p(t, y|x)-p(t, y|-x)
\nonumber\\
&=& \frac{1}{\sqrt{2 \pi t}}
\Big(
e^{-(y-x)^2/2t}-e^{-(y+x)^2/2t} \Big),
\quad x, y \in \R_+, \, t \geq 0.
\nonumber
\end{eqnarray}
If we put two absorbing walls at the origin and
at $x=h >0$, then
repeated application of the reflection principle
determines
the transition probability
density of the absorbing BM in the interval
$(0, h)$ as
\begin{eqnarray}
p_{2}^{h}(t, y|x)
&=& \sum_{n=-\infty}^{\infty} \Big\{
p(t, y|x+2hn)- p(t, y|-x+2hn) \Big\}
\nonumber\\
&=& \frac{1}{\sqrt{2 \pi t}}
\sum_{n=-\infty}^{\infty} \left[
\exp\left\{ -\frac{1}{2t}
(y-(x+2hn))^2 \right\}
- \exp\left\{ -\frac{1}{2t}
(y-(-x+2hn))^2 \right\}
\right] 
\nonumber
\end{eqnarray}
for $x, y \in (0, h), t \geq 0$.
Since \BES3, $X(t)$, is equivalent with the
BM conditioned to stay positive,
and this process is realized as an $h$-transform of
the absorbing BM with a wall at the origin
(see, for example, \cite{KT07b}),
we will see that
\begin{equation}
\P(H_1 < h)=\lim_{x \to 0, y \to 0}
\frac{p_{2}^{h}(1, y|x)}{p_1(1, y|x)}
\label{eqn:Ph0}
\end{equation}
for (\ref{eqn:H1}).
The limit of (\ref{eqn:Ph0}) can be readily performed
and we have
$$
\P(H_1 < h)=\sum_{n=-\infty}^{\infty}
e^{-2 h^2 n^2}(1-4 h^2 n^2)
$$
and the probability density is 
obtained as
\begin{eqnarray}
q_1(h) &\equiv& \frac{d}{dh} 
\P(H_1 < h) \nonumber\\
&=& 8 \sum_{n=1}^{\infty}
e^{-2 h^2 n^2}
(4 h^3 n^4-3 h n^2).
\nonumber
\end{eqnarray}
The $s$-th moment of $H_1$ is defined by
$$
\E[H_1^{s}]=\int_{0}^{\infty} dh \,
h^{s} q_1(h),
$$
and (\ref{eqn:equal1}) is derived, for which
the following equalities are useful,
\begin{eqnarray}
\int_{0}^{\infty} dh \, h^{s}
e^{-2h^2 n^2}
&=&2^{-(s+3)/2} n^{-(s+1)}
\Gamma((s+1)/2), \nonumber\\
\label{eqn:Gammaeq1}
\Gamma(s+1) &=& s \Gamma(s), 
\quad \Re s > 0.
\end{eqnarray}

\begin{figure}
\includegraphics[width=0.8\linewidth]{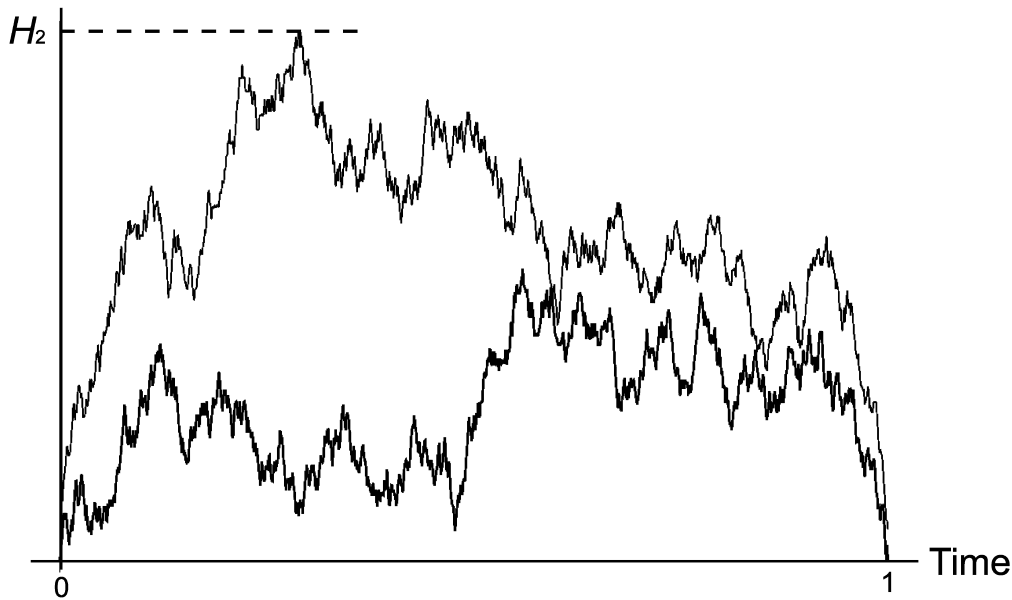}
\caption{Sample path of two Bessel bridges with duration 1
conditioned never to collide.
\label{fig:BridgeFig2}}
\end{figure}
For $N=2,3, \cdots$, the noncolliding $N$-particle 
system of Bessel bridges with duration 1,
with all particles started from 0, is denoted by
$\widetilde{\X}_N(t)
=(\widetilde{X}_{1}(t), \cdots, \widetilde{X}_N(t))$, 
where
$$
0 < \widetilde{X}_1(t) <
\widetilde{X}_2(t) < \cdots
< \widetilde{X}_N(t),
\quad 0 < t < 1.
$$
The stochastic variable $H_N$ is defined as the
maximum value of the $N$-th Bessel bridge
(\ref{eqn:HN}).
See Figure 2 for the $N=2$ case.
By the Karlin-McGregor formula \cite{KM59},
we will have
$$
\P(H_N < h)
=\lim_{x_j \to 0, y_j \to 0,
1 \leq j \leq N}
F_{h}(y_1, y_2, \cdots, y_N|
x_1, x_2, \cdots, x_N),
$$
where
\begin{equation}
F_{h}(y_1, y_2, \cdots, y_N|
x_1, x_2, \cdots, x_N)
\equiv \frac{\displaystyle{
\det_{1 \leq j, k \leq N}
[p_2^{h}(1, y_j | x_k)]}}
{\displaystyle{
\det_{1 \leq j, k \leq N}
[p_1(1, y_j| x_k)]}}
\label{eqn:Fh1}
\end{equation}
for $\x=(x_1, \cdots, x_N),
\y=(y_1, \cdots, y_N) 
\in \W_N^{h}
\equiv \{0 < x_1 < \cdots < x_N < h\}$.
The $s$-th moment of $H_N$ is then given by
$$
\E[H_N^{s}]=
\int_{0}^{\infty} dh \, 
h^s q_N(h) \quad
\mbox{with} \quad
q_N(h) \equiv \frac{d}{dh} \P(H_N < h).
$$

\subsection{Results}
Here we show our expressions for the moments
of $H_2$ of the two Bessel bridges 
conditioned never to collide in $t \in (0,1)$.
First expression is given using 
the double Dirichlet series of the form
\begin{equation}
Z(\alpha, \beta; \gamma)
\equiv 
\sum_{(n_1, n_2) \in \Z^2 \setminus \{(0,0)\}}
\frac{n_{1}^{\alpha} n_{2}^{\beta}}
{(n_{1}^2+n_{2}^{2})^{\gamma}},
\label{eqn:Z1}
\end{equation}
where $\Z$ denotes the set of all integers.
\begin{prop}
\label{thm:moment1}
Let 
\begin{equation}
\widetilde{Z}_{a}(b)=
\Gamma(a+2b) Z(2b,2b;a+2b).
\label{eqn:Ztil1}
\end{equation}
Then
\begin{eqnarray}
\E[H_2^s]
&=& \frac{2^{-s/2}}{24} s
\Big[ (s-1)(s^2-2s+12) \widetilde{Z}_{s/2}(0)
-4(s+4)(s+6) \widetilde{Z}_{s/2}(1)
+ 64 \widetilde{Z}_{s/2}(2) \Big].
\nonumber\\
\label{eqn:Pro2_1_eq}
\end{eqnarray}
\end{prop}
\vskip 0.5cm
\noindent{\bf Remark 1.} \,
As we show in Sect.3.1, 
$\widetilde{Z}_{s/2}(b), b=0,1,2$, have
simple poles at $s=0$ and $s=2$.
By the prefactor $s$ and by cancellation
of the $s=2$ poles among the three terms,
however, 
the expression (\ref{eqn:Pro2_1_eq}) has
finite limits in $s \to 0$ and $s \to 2$.
See Eq.(\ref{eqn:momentN2c}) below.
\vskip 0.5cm

We can rewrite this result using the Jacobi theta 
function (\ref{eqn:theta1}) and its derivatives, 
$
\vartheta'(u)=d \vartheta(u)/du,
\vartheta''(u)=d^2 \vartheta(u)/du^2.
$
\begin{thm}
\label{thm:main1}
Let 
\begin{equation}
K_{0}(s) =
\int_{1}^{\infty} du \, u^{s/2-1} 
\{\vartheta(u)^2-1 \},
\label{eqn:K0}
\end{equation}
and
\begin{eqnarray}
\xi_2(s) &=& -\frac{1}{6}
\left\{ (s+4)(s+6)
\int_{1}^{\infty} du \, u^{s/2+1}
\vartheta'(u)^2
\right.
\nonumber\\
&& \qquad \qquad \left.
+ ((2-s)+4)((2-s)+6)
\int_{1}^{\infty} du \, u^{(2-s)/2+1}
\vartheta'(u)^2
\right\}
\nonumber\\
&& + \frac{8}{3}
\int_{1}^{\infty} du \,
(u^{s/2+3}+u^{(2-s)/2+3})
\vartheta''(u)^2 
+\frac{1}{12}s(s-2) \vartheta(1)^2.
\label{eqn:xi2a}
\end{eqnarray}
Then
\begin{eqnarray}
\E[H_2^{s}]
&=& \left( \frac{\pi}{2} \right)^{s/2}
\left[ \frac{1}{24}(1-s)(s^2-2s+12)(2-s K_0(s))
\right. \nonumber\\
&& \left. \qquad \qquad 
-4 s \Big( \vartheta(1)
\vartheta'(1)
+2 s \vartheta'(1)^2
\Big) 
+s \xi_2(s) \right], \quad
s \in \C.
\label{eqn:momentN2d}
\end{eqnarray}
\end{thm}
\vskip 0.5cm
\noindent{\bf Remark 2.} \,
By the integral representation (\ref{eqn:xisym1}),
it is clear that $\xi(s)$ satisfies the functional
equation
$$
\xi(1-s)=\xi(s), \quad s \in \C.
$$
It is interesting to see that the function
$\xi_2(s)$, which appears in the expression
(\ref{eqn:momentN2d}),
satisfies the functional equation
$$
\xi_2(2-s)=\xi_2(s), \quad s \in \C.
$$
\vskip 0.5cm

As will be explicitly given in Sect.3.3,
$\xi(s), K_{0}(x)$ and $\xi_2(s)$
are expressed using series of 
the incomplete gamma functions.
Numerical evaluation of the incomplete
gamma functions is easy, and
the series converge rapidly.
Actually we have readily obtained the
values of moments for $N=1$
and $N=2$ as shown in Table 1.
(The trivial result $\E[1]=1$
is obtained by setting $s=0$ in (\ref{eqn:equal1})
with (\ref{eqn:xisym1}) and in (\ref{eqn:momentN2d}).
Since we know Euler's work on the relation
between $\zeta(2n), n=1,2,3, \cdots$ and
the Bernoulli numbers,
Eq.(\ref{eqn:equal1}) gives 
$\E[H_1^2]=\zeta(2)=\pi^2/6=1.644934 \cdots$ and
$\E[H_1^4]=3 \zeta(4)=\pi^4/30
=3.246969 \cdots$.)
By the relations (\ref{eqn:equal2}) 
and (\ref{eqn:limitA}) with $N=2$,
from the values in the $s=1$ column in the Table 1,
the dominant terms of the previous results
(\ref{eqn:BKR}) and (\ref{eqn:Fulmek1}) 
are reproduced;
\begin{eqnarray}
\langle h_1(2n) \rangle
&\simeq& \sqrt{2n} \times
\E[H_1] \nonumber\\
&=& \sqrt{2n} \times 1.253314 \cdots
= \sqrt{\pi n},
\nonumber\\
\langle h_2(2n) \rangle
&\simeq& \sqrt{2n} \times
\E[H_2] \nonumber\\
&=& \sqrt{2n} \times 1.822625 \cdots
= 2.57758 \cdots \times \sqrt{n}.
\nonumber
\end{eqnarray}
\begin{table}
\caption{Numerical values of moments}
\label{tab1:moment}
\begin{center}
\begin{tabular}{|c||r|r|r|r|r|r|}
\hline
$s$ & 0 & 1 & 2 & 3 & 4 & 5 \cr
\hline
\hline
$\E[H_1^s]$ & 1.0 & 1.253314 & 1.644934 & 2.259832 &
3.246969 & 4.873485 \cr
\hline
$\E[H_2^s]$ & 1.0 & 1.822625 & 3.395156 & 6.463823 &
12.576665 & 25.005999 \cr
\hline
\end{tabular}
\end{center}
\end{table}

\SSC{Proofs and Numerical Calculations}
\subsection{Proofs of Theorems}
In this subsection we will prove
Proposition \ref{thm:moment1} and 
Theorem \ref{thm:main1}.
Let $r_{\rm d}(y_1, y_2|x_1, x_2)$ and
$r_{\rm n}(y_1, y_2| x_1, x_2)$ be the denominator
and the numerator of the
RHS of (\ref{eqn:Fh1}), respectively,
given by $2 \times 2$ determinants,
when $N=2$.
We have obtained the estimations
\begin{eqnarray}
r_{\rm d}(x_1, x_2|x_1, x_2) &=& \frac{1}{3 \pi}
x_1^2 x_2^2 (x_1-x_2)^2
(x_1+x_2)^2
+O(x_1^{10}, x_2^{10}), \nonumber\\
r_{\rm n}(x_1, x_2|x_1, x_2) &=& \frac{1}{9 \pi}
x_1^2 x_2^2 (x_1-x_2)^2
(x_1+x_2)^2 
\sum_{n_1=-\infty}^{\infty} \sum_{n_2=-\infty}^{\infty}
e^{-2h^2(n_1^2+n_2^2)} Q_h(n_1, n_2)
\nonumber\\
&& \hskip 3cm+O(x_1^9, x_2^9), 
\nonumber
\end{eqnarray}
for $x_1, x_2 \ll 1$ with
\begin{eqnarray}
Q_h(n_1, n_2)
&=& 3-48h^2 n_1^2 + 72 h^4 n_1^4
+72 h^4 n_1^2 n_2^2
-32 h^6 n_1^6 -96 h^6 n_1^4 n_2^2
\nonumber\\
&& 
+128 h^8 n_1^6 n_2^2
-128 h^8 n_1^4 n_2^4.
\nonumber
\end{eqnarray}
Then the following results are concluded.
\begin{lem}
\label{thm:density}
$$
\P(H_2 < h) =
\sum_{(n_1, n_2) \in \Z^2}
e^{-2h^2(n_1^2+n_2^2)} A_h(n_1, n_2)
$$
with
\begin{eqnarray}
A_h(n_1, n_2)
&=& 1- 16 h^2 n_1^2+24 h^4 n_1^4
+ 24 h^4 n_1^2 n_2^2
-\frac{32}{3} h^6 n_1^6
- 32 h^6 n_1^4 n_2^2
\nonumber\\
&& +\frac{128}{3} h^8 n_1^6 n_2^2
-\frac{128}{3} h^8 n_1^4 n_2^4.
\nonumber
\end{eqnarray}
And then
$$
q_2(h)= 
\sum_{(n_1, n_2) \in \Z^2 \setminus \{(0,0)\}}
e^{-2 h^2(n_1^2+n_2^2)}
B_h(n_{1}, n_{2})
$$
with
\begin{eqnarray}
B_h(n_1, n_2)
&=& \frac{8}{3} h \Big\{
-15 n_1^2 + 60 h^2 n_1^4 + 60 h^2 n_1^2 n_2^2
-60 h^4 n_1^6 -180 h^4 n_1^4 n_2^2
\nonumber\\
&& +16 h^6 n_1^8 + 192 h^6 n_1^6 n_2^2
-80 h^6 n_1^4 n_2^4
-64 h^8 n_1^8 n_2^2 +64 h^8 n_1^6 n_2^4 \Big\}.
\nonumber
\end{eqnarray}
\end{lem}
\vskip 0.5cm

\noindent{\it Proof of Proposition \ref{thm:moment1}}. \,
In order to describe the moments
$\E[H_2^s]=\int_0^{\infty} dh \, h^s q_2(h)$, 
we introduce the following notation,
\begin{equation}
I_{s}(\alpha, \beta)
=\sum_{(n_1, n_2) \in \Z^2 \setminus \{(0,0)\}}
n_{1}^{\alpha} n_{2}^{\beta}
\int_{0}^{\infty} dh \, h^{\alpha+\beta-1+s}
e^{-2h^2(n_{1}^2+n_{2}^2)}.
\label{eqn:Iab1}
\end{equation}
Then, by Lemma \ref{thm:density}, we have
\begin{eqnarray}
\E[H_2^s] &=& \frac{8}{3} \Big\{
-15 I(2,0)+60 I(4,0) +60 I(2,2) -60 I(6,0)
-180 I(4,2) \nonumber\\
&& + 16 I(8,0) + 192 I(6,2) -80 I(4,4)
-64 I(8,2) + 64 I(6,4) \Big\}.
\label{eqn:moment2}
\end{eqnarray}
Next we rewrite $I_s(\alpha, \beta)$ using
the Gamma function (\ref{eqn:gamma1}).
In the integrals in (\ref{eqn:Iab1}), 
we change the integral variable from $h$ to $u$ by
$
u=2h^2(n_1^2+n_2^2), 
$
respectively.
Then we have
$$
I_{s}(\alpha, \beta)
= 2^{-(\alpha+\beta+2+s)/2}
\Gamma((\alpha+\beta+s)/2)
Z( \alpha, \beta; (\alpha+\beta+s)/2),
$$
where $Z(\alpha, \beta; \gamma)$ is defined by (\ref{eqn:Z1}).
From (\ref{eqn:moment2}) we will see that
\begin{eqnarray}
\E[H_2^s] &=& \frac{1}{3} 2^{(2-s)/2}
\Big[ -15 \Gamma(1+s/2) Z(2, 0; 1+s/2)
\nonumber\\
&&
+30 \Gamma(2+s/2) 
\{ Z(4,0; 2+s/2)+Z(2,2; 2+s/2)\}
\nonumber\\
&& -15 \Gamma(3+s/2) 
\{ Z(6,0; 3+s/2) + 3 Z(4,2; 3+s/2) \}
\nonumber\\
&&
+2 \Gamma(4+s/2) 
\{ Z(8,0; 4+s/2) + 12 Z(6,2; 4+s/2)
-5 Z(4,4; 4+s/2)\}
\nonumber\\
&& -4 \Gamma(5+s/2)
\{ Z(8,2; 5+s/2)
-Z(6,4; 5+s/2) \} \Big].
\label{eqn:moment3}
\end{eqnarray}
By definition (\ref{eqn:Z1}), 
\begin{eqnarray}
&& Z(2,0; 1+s/2)
= \frac{1}{2} \sum_{(n_1, n_2) \in \Z^2 \setminus \{(0,0)\}}
\frac{n_{1}^{2}+n_{2}^{2}}
{(n_{1}^2+n_{2}^{2})^{1+s/2}}
= \frac{1}{2} Z(0,0; s/2),
\nonumber\\
&& Z(4,0;2+s/2)+Z(2,2; 2+s/2)
= \frac{1}{2} 
\sum_{(n_1, n_2) \in \Z^2 \setminus \{(0,0)\}}
\frac{n_1^4+2n_1^2 n_2^2+n_2^4}
{(n_{1}^2+n_{2}^{2})^{2+s/2}}
\nonumber\\
&& \qquad \qquad = \frac{1}{2} Z(0,0; s/2),
\nonumber
\end{eqnarray}
\begin{eqnarray}
&& Z(6,0;3+s/2)+3Z(4,2; 3+s/2)
= \frac{1}{2} 
\sum_{(n_1, n_2) \in \Z^2 \setminus \{(0,0)\}}
\frac{n_1^6+3n_1^4 n_2^2+3 n_1^2 n_2^4 + n_2^6}
{(n_{1}^2+n_{2}^{2})^{3+s/2}}
\nonumber\\
&& \qquad \qquad =\frac{1}{2} Z(0,0;s/2),
\nonumber
\end{eqnarray}
\begin{eqnarray}
&& Z(8,0; 4+s/2)+12 Z(6,2; 4+s/2)
-5 Z(4,4; 4+s/2) \nonumber\\
&& \qquad \qquad
= \sum_{(n_1, n_2) \in \Z^2 \setminus \{(0,0)\}}
\frac{n_{1}^8+12 n_1^6 n_2^2 - 5 n_1^4 n_2^4}
{(n_{1}^2+n_{2}^{2})^{4+s/2}}
\nonumber\\
&& \qquad \qquad = \frac{1}{2}
\sum_{(n_1, n_2) \in \Z^2 \setminus \{(0,0)\}}
\frac{(n_1^2+n_2^2)^4+
8 n_1^2 n_2^2 (n_1^2+n_2^2)^2
-32 n_1^4 n_2^4}
{(n_{1}^2+n_{2}^{2})^{4+s/2}}
\nonumber\\
&& \qquad \qquad = \frac{1}{2} Z(0,0; s/2)
+ 4 Z(2,2; 2+s/2)
-16 Z(4,4; 4+s/2),
\nonumber
\end{eqnarray}
\begin{eqnarray}
&&Z(8,2; 5+s/2)-Z(6,4; 5+s/2) 
=\sum_{(n_1, n_2) \in \Z^2 \setminus \{(0,0)\}}
\frac{n_{1}^8 n_2^2-n_1^6 n_2^4}
{(n_{1}^2+n_{2}^{2})^{5+s/2}}
\nonumber\\
&& \qquad \qquad = \frac{1}{2}
\sum_{(n_1, n_2) \in \Z^2 \setminus \{(0,0)\}}
\frac{n_1^2 n_2^2(n_1^2+n_2^2)^3
-4 n_1^4 n_2^4 (n_1^2+n_2^2)
}
{(n_{1}^2+n_{2}^{2})^{5+s/2}}
\nonumber\\
&& \qquad \qquad = \frac{1}{2} Z(2,2; 2+s/2)
-2 Z(4,4; 4+s/2).
\nonumber
\end{eqnarray}
Using the above results, 
(\ref{eqn:moment3}) is rewritten as
$$
\E[H_2^s]
= \frac{2^{(2-s)/2}}{3}
\Big[ c_1(s) Z(0,0; s/2)
+ c_2(s) Z(2,2; 2+s/2)
+ c_3(s) Z(4,4; 4+s/2) \Big]
$$
with
\begin{eqnarray}
c_1(s) &=& -\frac{15}{2} \Gamma(1+s/2)
+15 \Gamma(2+s/2)
-\frac{15}{2} \Gamma(3+s/2)
+\Gamma(4+s/2)
\nonumber\\
&=& \frac{1}{16} s(s-1) (s^2-2s+12) \Gamma(s/2),
\nonumber\\
c_2(s) &=& 8 \Gamma(4+s/2)-2 \Gamma(5+s/2)
\nonumber\\
&=& -s \Gamma(s/2+4) \nonumber\\
&=& - \frac{1}{4} s(s+4)(s+6) \Gamma(s/2+2),
\nonumber\\
c_3(s) &=& -32 \Gamma(4+s/2)+ 8\Gamma(5+s/2)
\nonumber\\
&=& 4s \Gamma(s/2+4),
\nonumber
\end{eqnarray}
where (\ref{eqn:Gammaeq1}) has been used.
Then Proposition \ref{thm:moment1} was proved.
\qed 
\vskip 0.5cm
Let ${\bf 1}_{\{\omega\}}$ be the
indicator function of condition $\omega$;
${\bf 1}_{\{\omega\}}=1$, if
$\omega$ is satisfied and
${\bf 1}_{\{\omega\}}=0$, otherwise.
The following equality is derived.
\begin{lem}
\label{thm:Z_theta}
\begin{equation}
\widetilde{Z}_a(b)
=\pi^{a} 
\int_{0}^{\infty} du \,
u^{a+2b-1} \left\{
\left( \frac{d^b}{du^b} 
\vartheta(u) \right)^2
-{\bf 1}_{\{b=0\}} \right\}.
\label{eqn:Ztil2}
\end{equation}
\end{lem}
\noindent{\it Proof.} \,
By definition (\ref{eqn:Ztil1})
\begin{equation}
\widetilde{Z}_{a}(b)
= \sum_{(n_1, n_2) \in \Z^2 \setminus \{(0,0)\}}
\frac{n_1^{2b} n_2^{2b}}{(n_1^2+n_2^2)^{a+2b}}
\Gamma(a+2b).
\label{eqn:Ztil3}
\end{equation}
By changing the integral variable $u$ by $w$
with $u=\pi(n_1^2+n_2^2)w$ in the integral
(\ref{eqn:gamma1}), we have
$$
\Gamma(s)=\pi^{s} (n_1^2+n_2^2)^s
\int_{0}^{\infty} dw \, w^{s-1}
e^{-\pi(n_1^2+n_2^2)w}.
$$
Then (\ref{eqn:Ztil3}) becomes
\begin{eqnarray}
\widetilde{Z}_{a}(b)
&=& \pi^{a+2b} \int_{0}^{\infty} dw \,
w^{a+2b-1} 
\sum_{(n_1,n_2) \in \Z^2 \setminus \{(0,0)\}}
n_1^{2b} n_2^{2b} e^{-\pi(n_1^2+n_2^2)w}
\nonumber\\
&=& \pi^{a+2b} \int_{0}^{\infty} dw \,
w^{a+2b-1}
\left\{ \left( \sum_{n=-\infty}^{\infty} n^{2b}
e^{-\pi n^2 w} \right)^2 
-{\bf 1}_{\{b=0\}} \right\}.
\nonumber
\end{eqnarray}
Since $\displaystyle{(-\pi n^2)^{b}e^{-\pi n^2 w}
=\frac{d^b}{dw^b} e^{-\pi n^2 w}}$,
(\ref{eqn:Ztil2}) is obtained. \qed
\vskip 0.5cm
Then we have the following expressions of 
moments.
\begin{prop}
\label{thm:moment2}
\begin{eqnarray}
\E[H_2^s]
&=& \frac{1}{24} \left(\frac{\pi}{2}\right)^{s/2}
s \left[ (s-1)(s^2-2s+12)
\int_{0}^{\infty} du \, u^{s/2-1}
\{\vartheta(u)^2 -1 \}
\right.
\nonumber\\
&& \left.
-4(s+4)(s+6)
\int_{0}^{\infty} du \,
u^{s/2+1} \vartheta'(u)^2
+ 64 \int_{0}^{\infty} du \,
u^{s/2+3} \vartheta''(u)^2
\right].
\nonumber\\
\label{eqn:momentN2b}
\end{eqnarray}
\end{prop}
\vskip 0.5cm
\noindent{\it Proof of Theorem \ref{thm:main1}.} \,
By the reciprocity law of the Jacobi theta function
\cite{BPY01}
\begin{equation}
\vartheta(u)
=\sqrt{\frac{1}{u}}
\vartheta\left(\frac{1}{u}\right),
\quad \Re u >0,
\label{eqn:recip1}
\end{equation}
we can show the following, 
\begin{eqnarray}
I_{1} &\equiv& \int_{0}^{\infty} du \, u^{s/2-1} 
\{ \vartheta(u)^2 -1 \}
\nonumber\\
&=& -\frac{2}{s}+\frac{2}{s-2}
+\int_{1}^{\infty} du \, u^{-s/2} 
\{ \vartheta(u)^2 -1 \}
+\int_{1}^{\infty} du \, u^{s/2-1} 
\{ \vartheta(u)^2 -1 \},
\label{eqn:I1A}
\end{eqnarray}
\begin{eqnarray}
I_{2} &\equiv& \int_{0}^{\infty} du \,
u^{s/2+1} \vartheta'(u)^2
\nonumber\\
&=& \frac{1}{2(s-2)} +\int_{1}^{\infty} du \, u^{s/2+1}
\vartheta'(u)^2
+\int_{1}^{\infty} du \, u^{-s/2+1}
\vartheta(u) \vartheta'(u)
\nonumber\\
&& +\int_{1}^{\infty} du \, u^{-s/2+2}
\vartheta'(u)^2
+\frac{1}{4} \int_{1}^{\infty} du \, u^{-s/2}
\{\vartheta(u)^2-1 \},
\label{eqn:I2A}
\end{eqnarray}
\begin{eqnarray}
I_{3} &\equiv& \int_{0}^{\infty} du \,
u^{s/2+3} \vartheta''(u)^2
\nonumber\\
&=& \frac{9}{8(s-2)}
+\int_{1}^{\infty} du \, u^{s/2+3}
\vartheta''(u)^2
+\int_{1}^{\infty} du \, u^{-s/2+4}
\vartheta''(u)^2
\nonumber\\
&& +6 \int_{1}^{\infty} du \, u^{-s/2+3}
\vartheta'(u) \vartheta''(u)
+\frac{3}{2} \int_{1}^{\infty} du \, u^{-s/2+2}
\vartheta(u) \vartheta''(u)
\nonumber\\
&& +\frac{9}{2} \int_{1}^{\infty} du \, u^{-s/2+1}
\vartheta(u) \vartheta'(u) 
+9 \int_{1}^{\infty} du \, u^{-s/2+2}
\vartheta'(u)^2
\nonumber\\
&& 
+\frac{9}{16} \int_{1}^{\infty} du \, u^{-s/2}
\{ \vartheta(u)^2-1 \}.
\label{eqn:I3A}
\end{eqnarray}
The derivations are given in Appendix A.
They show that $I_1$ has only two simple poles at 
$s=0$ (with residue $-2$) and at $s=2$ (with residue $2$),
$I_2$ does only one simple pole at
$s=2$ (with residue $1/2$),
and $I_3$ does only one simple pole at
$s=2$ (with residue $9/8$).
Using them in the representation (\ref{eqn:momentN2b}),
we have 
\begin{eqnarray}
\E[ H_2^{s} ]
&=& \frac{1}{24} \left( \frac{\pi}{2} \right)^{s/2}
\Bigg[ 2(s^2-14s+12) 
\nonumber\\
&& \quad
+s \Big\{ (s-1)(s^2-2s+12) K_{0}(s)
-4(s+4)(s+6) K_{1}(s) + 64 K_{2}(s)
\nonumber\\
&& \qquad +s(s-2)^2 K_{0}(2-s)
-4(s^2+10s-120) K_{1}(2-s) 
+64 K_{2}(2-s)
\nonumber\\
&& \qquad -4(s^2+10s-48) J_{1}(s) +96 J_{2}(s)
+384 J_{3}(s)
\Big\} \Bigg],
\label{eqn:momentN2c}
\end{eqnarray}
where
\begin{eqnarray}
K_{0}(s) &=&
\int_{1}^{\infty} du \, u^{s/2-1} 
\{\vartheta(u)^2-1 \},
\nonumber\\
K_{1}(s) &=&
\int_{1}^{\infty} du \, u^{s/2+1}
\vartheta'(u)^2,
\nonumber\\
K_{2}(s) &=&
\int_{1}^{\infty} du \, u^{s/2+3}
\vartheta''(u)^2,
\nonumber\\
J_{1}(s) &=& 
\int_{1}^{\infty} du \, u^{-s/2+1}
\vartheta(u) \vartheta'(u),
\nonumber\\
J_{2}(s) &=& \int_{1}^{\infty} du \, u^{-s/2+2}
\vartheta(u) \vartheta''(u)
\nonumber\\
J_{3}(s) &=&
\int_{1}^{\infty} du \, u^{-s/2+3}
\vartheta'(u)
\vartheta''(u).
\nonumber
\end{eqnarray}
Note that the singularities in $I_1, I_2, I_3$ have 
been all cancelled 
and all integrals converge for all $s \in \C$.

Now we perform partial integrations
to evaluate $J_1(s), J_2(s)$ and $J_3(s)$.
For $J_1(s)$,
\begin{eqnarray}
J_1(s) &=& \Big[ u^{-s/2+1}
\vartheta(u)^2 \Big]_{1}^{\infty}
-\int_{1}^{\infty} du \,
\frac{d}{du} \Big(
u^{-s/2+1} \vartheta(u) \Big)
\vartheta(u)
\nonumber\\
&=& \Big[ u^{-s/2+1}
\vartheta(u)^2 \Big]_{1}^{\infty}
+\frac{1}{2}(s-2) K_0(2-s)
-\Big[ u^{-s/2+1} \Big]_{1}^{\infty}- J_1(s).
\nonumber
\end{eqnarray}
Since
$$
\Big[ u^{-s/2+1}
\vartheta(u)^2 \Big]_{1}^{\infty}
-\Big[ u^{-s/2+1} \Big]_{1}^{\infty}
= -\{\vartheta(1)^2-1 \},
$$
we have
\begin{equation}
J_{1}(s)=\frac{1}{4} (s-2) K_0(2-s)
-\frac{1}{2} \{\vartheta(1)^2-1 \}.
\label{eqn:J1}
\end{equation}
For $J_2(s)$, we see
\begin{eqnarray}
J_{2}(s) &=& \left[ u^{-s/2+2} \vartheta(u)
\frac{d}{du} \vartheta(u) \right]_{1}^{\infty}
- \int_{1}^{\infty} du \,
\frac{d}{du} \Big\{
u^{-s/2+2} \vartheta(u) \Big\}
\frac{d}{du} \vartheta(u)
\nonumber\\
&=& - \vartheta(1) \vartheta'(1)
+ \frac{1}{2}(s-4)J_1(s) - K_1(2-s).
\nonumber
\end{eqnarray}
Inserting (\ref{eqn:J1}), we have
\begin{equation}
J_2(s) =
- \vartheta(1) \vartheta'(1)
- \frac{1}{4}(s-4)
\{\vartheta(1)^2-1 \}
+\frac{1}{8}(s-2)(s-4) K_0(2-s) -K_1(2-s).
\label{eqn:J2}
\end{equation}
Similarly by partial integration, we obtain
\begin{eqnarray}
J_3(s) &=&
\left[ u^{-s/2+3} 
\left(\frac{d}{du} \vartheta(u)\right)^2
\right]_{1}^{\infty}
+\frac{1}{2}(s-6) K_1(2-s) -J_3(s)
\nonumber\\
&=& -\frac{1}{2} \vartheta'(1)^2
+\frac{1}{4}(s-6) K_1(2-s).
\label{eqn:J3}
\end{eqnarray}
Using (\ref{eqn:J1})-(\ref{eqn:J3})
in (\ref{eqn:momentN2c}), 
(\ref{eqn:momentN2d}) with
(\ref{eqn:K0}) and (\ref{eqn:xi2a}) is obtained.
\qed

\subsection{Expressions by Incomplete Gamma Functions}

Let $\Gamma(z,p)$ be the incomplete gamma function of
the second kind defined as
\begin{eqnarray}
\Gamma(z, p) &=& \int_{p}^{\infty} du \, 
u^{z-1} e^{-u}
\nonumber\\
&=& \Gamma(z)-\int_{0}^{p} du \, u^{z-1} e^{-u},
\quad \Re z >0, \, p > 0.
\label{eqn:ingamma}
\end{eqnarray}
As demonstrated in Appendix B, the integrals
appearing in (\ref{eqn:xisym1}) and in our result
given by Theorem \ref{thm:main1} are expressed using
$\Gamma(z,p)$.
We have obtained the following expressions.
\begin{eqnarray}
\E[H_1^s]
&=& \left( \frac{\pi}{2} \right)^{s/2}
\Bigg[ 1+ s(s-1) \nonumber\\
&& \times \left. \left\{ 
\pi^{-s/2}\sum _{n=1}^{\infty }n^{-s}
 \Gamma(s/2,\pi n^2)
 + \pi^{s/2-1/2}\sum _{n=1}^{\infty }n^{s-1}
 \Gamma(-s/2+1/2,\pi n^2) \right\}
\right], \qquad
\label{eqn:N1momentGamma}
\end{eqnarray}
and
\begin{eqnarray}
\E[H_2^s] 
 &=& \left( \frac{\pi}{2} \right)^{s/2}
 \Bigg[ \frac{1}{12}(1-s)(s^2-2s+12) \nonumber\\
&& 
 \times \left\{
 1-2 s\pi^{-s/2} \left( 
  \sum_{n_1=1}^{\infty} \sum_{n_2=1}^{\infty}
\frac{1}{(n_1^2+n_2^2)^{s/2}}
\Gamma( s/2, \pi (n_1^2+n_{2}^2))
+ \sum_{n=1}^{\infty}
\frac{1}{n^s} \Gamma(s/2, \pi n^2)
 \right)
  \right\} \nonumber\\
&& \quad 
-4 s \{ \vartheta(1)
\vartheta'(1)
+2 (\vartheta'(1))^2 \} 
+s \xi_2(s) 
\Bigg],
\label{eqn:N2momentGamma}
\end{eqnarray}
with
\begin{eqnarray}
\xi_2(s)
 &=& -\frac{2}{3}
\left\{ \pi^{-s/2}(s+4)(s+6)
\sum_{n_1=1}^{\infty} \sum_{n_2=1}^{\infty}
\frac{n_1^2 n_2^2}
{(n_1^2+n_2^2)^{s/2+2}}
\Gamma( s/2+2, \pi (n_1^2+n_2^2))
 \right. \nonumber\\
&& \left.
 +\pi^{s/2-1}((2-s)+4)((2-s)+6)
 \sum_{n_1=1}^{\infty} \sum_{n_2=1}^{\infty}
\frac{n_1^2 n_2^2}
{(n_1^2+n_2^2)^{-s/2+3}}
\Gamma( -s/2+3, \pi (n_1^2+n_2^2))
\right\}
\nonumber\\
&&
 +\frac{32}{3}\left\{
 \pi^{-s/2}
\sum_{n_1=1}^{\infty} \sum_{n_2=1}^{\infty}
\frac{n_1^4 n_2^4}
{(n_1^2+n_2^2)^{s/2+4}}
\Gamma( s/2+4, \pi (n_1^2+n_2^2))
\right. \nonumber\\
&& \left. \qquad 
+ \pi^{s/2-1}
\sum_{n_1=1}^{\infty} \sum_{n_2=1}^{\infty}
\frac{n_1^4 n_2^4}
{(n_1^2+n_2^2)^{-s/2+5}}
\Gamma( -s/2+5, \pi (n_1^2+n_2^2))
  \right\}
 \nonumber\\
 && +\frac{1}{12}s(s-2) \vartheta(1)^2.
 \label{eqn:FsGamma}
\end{eqnarray}
The values in Table 1 were obtained
by numerically calculating the
incomplete gamma functions and the series
of them with appropriate coefficients
in the above expressions.

\vskip 1cm
\appendix
\begin{LARGE}
{\bf Appendices}
\end{LARGE}
\SSC{Calculation of $I_j, j=1,2,3$}
By the reciprocity law of the Jacobi theta function 
(\ref{eqn:recip1}),
\begin{eqnarray}
I_{1} &=& 
\int_{0}^{1} du \, u^{s/2-1}
\{\vartheta(u)^2-1 \}
+\int_{1}^{\infty} du \, u^{s/2-1}
\{\vartheta(u)^2-1 \}
\nonumber\\
&=& \int_{0}^{1} du \, u^{s/2-1}
\left\{ \frac{1}{u} \left(
\vartheta \left(\frac{1}{u}\right) \right)^2 
-1 \right\}
+\int_{1}^{\infty} du \, u^{s/2-1}
\{\vartheta(u)^2-1 \}
\nonumber\\
&=& \int_{0}^{1} du \, u^{s/2-1}
\left\{ \frac{1}{u}-1 \right\}
\nonumber\\
&& +\int_{0}^{1} du \, u^{s/2-2} \left\{ \left(
\vartheta \left(\frac{1}{u}\right) \right)^2 -1 \right\}
+\int_{1}^{\infty} du \, u^{s/2-1}
\{\vartheta(u)^2-1 \}.
\nonumber
\end{eqnarray}
By calculating the first integral and by
changing the integral variable $u$ by $w=1/u$ in the
second integral, we obtain (\ref{eqn:I1A}).

Set
$$
I_{2}=\int_{0}^{1} du \, u^{s/2+1}
\left( \frac{d}{du} \vartheta(u) \right)^2
+ \int_{1}^{\infty} du \, u^{s/2+1}
\left( \frac{d}{du} \vartheta(u) \right)^2.
$$
In the first integrand, we use the reciprocity
law (\ref{eqn:recip1}) as
\begin{eqnarray}
&& \left( \frac{d}{du} \vartheta(u) \right)^2
= \left\{ \frac{d}{du} \left(
\sqrt{\frac{1}{u}} \vartheta
\left(\frac{1}{u} \right) \right) \right\}^2
\nonumber\\
&& \quad = 
\left\{ - \frac{1}{2} u^{-3/2} 
\vartheta
\left( \frac{1}{u} \right)
+\sqrt{\frac{1}{u}} \frac{d}{du}
\vartheta \left(\frac{1}{u}\right)
\right\}^2
\nonumber\\
&& \quad
= \frac{1}{4}u^{-3} \left(
\vartheta \left(\frac{1}{u}\right) \right)^2
-u^{-2} \vartheta\left(\frac{1}{u}\right)
\frac{d}{du} \vartheta\left(\frac{1}{u}\right)
+ u^{-1}
\left( \frac{d}{du} \vartheta\left(\frac{1}{u}\right)
\right)^2.
\nonumber
\end{eqnarray}
By inserting this into the first integral in $I_2$
and by changing the integral variables
$u$ by $w=1/u$, we have
\begin{eqnarray}
I_{2} &=& \frac{1}{4} \int_{1}^{\infty}
dw \, w^{-s/2} \vartheta(w)^2
+ \int_{1}^{\infty} dw \, w^{-s/2+1}
\vartheta(w) \frac{d}{dw} \vartheta(w)
\nonumber\\
&& + \int_{1}^{\infty} dw \, w^{-s/2+2}
\left( \frac{d}{dw} \vartheta(w) \right)^2
+ \int_{1}^{\infty} dw \, w^{s/2+1}
\left( \frac{d}{dw} \vartheta(w) \right)^2.
\nonumber
\end{eqnarray}
The first integral is rewritten as
\begin{eqnarray}
&& \frac{1}{4} \int_{1}^{\infty}
dw \, w^{-s/2} \vartheta(w)^2
= \frac{1}{4} \int_{1}^{\infty} dw \, w^{-s/2}
+\frac{1}{4} \int_{1}^{\infty}
dw \, w^{-s/2} \{
\vartheta(w)^2-1 \}
\nonumber\\
&& \quad =\frac{1}{2(s-2)}
+\frac{1}{4} \int_{1}^{\infty}
dw \, w^{-s/2} \{
\vartheta(w)^2-1 \}.
\nonumber
\end{eqnarray}
Then we obtain (\ref{eqn:I2A}).
For $I_{3}$ we set
$$
I_{3}=\int_{0}^{1} du \, u^{s/2+3} 
\left( \frac{d^2}{du^2} \vartheta(u) \right)^2
+ \int_{1}^{\infty} du \, u^{s/2+3} 
\left( \frac{d^2}{du^2} \vartheta(u) \right)^2
$$
and the first integrand is written using
the reciprocity law (\ref{eqn:recip1}) as
\begin{eqnarray}
&& 
\left( \frac{d^2}{du^2} \vartheta(u) \right)^2
=\left\{ \frac{d^2}{du^2}
\left( \sqrt{\frac{1}{u}}
\vartheta\left( \frac{1}{u} \right) \right) \right\}^2
\nonumber\\
&& \quad =
\left\{ \frac{3}{4} u^{-5/2}
\vartheta\left(\frac{1}{u}\right)
-u^{-3/2} \frac{d}{du} \vartheta
\left( \frac{1}{u} \right)
+u^{-1/2} \frac{d^2}{du^2}
\vartheta \left(\frac{1}{u}\right) \right\}^2
\nonumber\\
&& \quad =
\frac{9}{16}u^{-5} \left( \vartheta
\left(\frac{1}{u}\right) \right)^2
+u^{-3}\left(\frac{d}{du} \vartheta
\left(\frac{1}{u}\right) \right)^2
+u^{-1}\left(\frac{d^2}{du^2} \vartheta
\left(\frac{1}{u}\right) \right)^2
\nonumber\\
&& \quad -\frac{3}{2} u^{-4} 
\vartheta\left(\frac{1}{u}\right)
\frac{d}{du} \vartheta
\left(\frac{1}{u}\right)
+\frac{3}{2} u^{-3} 
\vartheta\left(\frac{1}{u}\right)
\frac{d^2}{du^2} \vartheta
\left(\frac{1}{u}\right)
\nonumber\\
&& \quad
-2 u^{-2} 
\left( \frac{d}{du}
\vartheta\left(\frac{1}{u}\right) \right)
\left(
\frac{d^2}{du^2} \vartheta
\left(\frac{1}{u}\right) \right).
\nonumber
\end{eqnarray}
By similar calculation we can obtain (\ref{eqn:I3A}).

\SSC{Incomplete Gamma Functions}

As an example, here we only
consider the integral
$$
L=\int_{1}^{\infty} du \, u^{-1/2}
\{\vartheta(u)^2-1 \}.
$$
By definition of the Jacobi theta function
(\ref{eqn:theta1}), 
we see that
\begin{eqnarray}
L &=& \int_{1}^{\infty} du \, u^{-1/2}
\sum_{(n_1, n_2) \in \Z^2 \setminus \{(0,0)\}}
e^{-\pi (n_1^2+n_2^2) u}
\nonumber\\
&=& 4 \int_{1}^{\infty} du \, u^{-1/2}
\sum_{n_1=1}^{\infty} \sum_{n_2=1}^{\infty}
e^{-\pi(n_1^2+n_2^2) u}
+ 4 \int_{1}^{\infty} du \, u^{-1/2}
\sum_{n=1}^{\infty} e^{-\pi n^2 u}.
\nonumber
\end{eqnarray}
In the first and the second integrals in the above,
we set 
$\pi(n_1^2+n_2^2)u =w$ and
$\pi n^2 u=w$, respectively, and replace the integral variables
$u$ by $w$. Then we have
\begin{eqnarray}
L &=& 4 \pi^{-1/2} \sum_{n_1=1}^{\infty} \sum_{n_2=1}^{\infty}
\frac{1}{(n_1^2+n_2^2)^{1/2}}
\int_{\pi(n_1^2+n_2^2)}^{\infty} dw
\, w^{-1/2} e^{-w}
+ 4 \pi^{-1/2} \sum_{n=1}^{\infty} \frac{1}{n}
\int_{\pi n^2}^{\infty} dw \,
w^{-1/2} e^{-w}
\nonumber\\
&=& 4 \pi^{-1/2} \sum_{n_1=1}^{\infty} \sum_{n_2=1}^{\infty}
\frac{1}{(n_1^2+n_2^2)^{1/2}}
\Gamma(1/2, \pi (n_1^2+n_{2}^2))
+ 4 \pi^{-1/2} \sum_{n=1}^{\infty}
\frac{1}{n} \Gamma(1/2, \pi n^2),
\nonumber
\end{eqnarray}
where $\Gamma(z,p)$ is the incomplete
gamma function defined by (\ref{eqn:ingamma}).
Combination of
similar calculations will give the expressions
(\ref{eqn:N1momentGamma})-(\ref{eqn:FsGamma}).

\vskip 1cm
\noindent{\bf Acknowledgements} \,
The authors would like to thank
Thomas Feierl for valuable comments
on this paper and 
on the recent progress of combinatorial methods
on the present problem.
M.K. is supported in part by
the Grant-in-Aid for Scientific Research 
(KIBAN-C, No.17540363) of Japan Society for
the Promotion of Science.

\end{document}